\documentclass{amsart}
\usepackage[utf8]{inputenc}
\usepackage[english]{babel}
\usepackage{amsmath}
\usepackage{amsthm}
\usepackage{tikz}
\usetikzlibrary{svg.path}
\usetikzlibrary{shapes,arrows,chains}
\usepackage{amsfonts}
\usepackage{amssymb}
\usepackage{graphicx}
\usepackage{xcolor}
\usepackage{bbm}
\usepackage{csquotes}
\usepackage{latexsym}
\usepackage{float}
\usepackage[backend=biber, style=alphabetic, sorting=nyt]{biblatex}
\newtheorem{thm}{Theorem}[section]
\newtheorem*{thm*}{Theorem}
\newtheorem{lem}[thm]{Lemma}

\newtheorem{prop}[thm]{Proposition}

\newtheorem{conj}[thm]{Conjecture}

\theoremstyle{definition}
\newtheorem{defn}[thm]{Definition}

\newtheorem{rem}[thm]{Remark}

\newtheorem{basass}[thm]{Basic Assumptions}

\newcommand{\codim}{\textnormal{codim}}

\renewcommand{\exp}{\textnormal{exp}}

\renewcommand{\o}[1]{\overline{ #1 }}

\newcommand{\im}{\textnormal{im}}

\renewcommand{\epsilon}{\varepsilon}
\renewcommand{\phi}{\varphi}

\newcommand{\sqvecttwo}[3]{\left[ \begin{matrix} #1 \\ #2 \\ #3 \end{matrix} \right]}
\newcommand{\sqvecteight}[9]{\left[ \begin{matrix} #1 \\ #2 \\ #3 \\#4 \\ #5 \\ #6 \\ #7 \\ #8 \\ #9 \end{matrix} \right]}

\newcommand{\vecttwo}[2]{\left(\begin{matrix} #1 \\ #2 \end{matrix} \right)}

\newcommand{\ctimes}{\mathbb{C}^\times}

\renewcommand{\Re}{\textnormal{Re}}
\renewcommand{\Im}{\textnormal{Im}}

\newcommand{\dhaus}{d_{\textnormal{Haus}}}

\newcommand{\hol}{\textnormal{hol}}
\title{On Some Systems of Equations in Abelian Varieties}
\author{Francesco Paolo Gallinaro}
\thanks{Address: Mathematisches Institut, Albert-Ludwigs-Universit\"at Freiburg, Ernst-Zermelo-Str. 1, 79100 Freiburg, Germany. Email address: francesco.gallinaro@mathematik.uni-freiburg.de}
\begin{document}
	
\maketitle

\begin{abstract}
We solve a case of the \textit{Abelian Exponential-Algebraic Closedness Conjecture}, a conjecture due to Bays and Kirby, building on work of Zilber, which predicts sufficient conditions for systems of equations involving algebraic operations and the exponential map of an abelian variety to be solvable in the complex numbers. More precisely, we show that the conjecture holds for subvarieties of the tangent bundle of an abelian variety $A$ which split as the product of a linear subspace of the Lie algebra of $A$ and an algebraic variety. This is motivated by work of Zilber and of Bays-Kirby, which establishes that a positive answer to the conjecture would imply quasiminimality of certain structures on the complex numbers. Our proofs use various techniques from homology (duality between cup product and intersection), differential topology (transversality) and o-minimality (definability of Hausdorff limits), hence we have tried to give a self-contained exposition.
\end{abstract}

\section{Introduction}

\let\thefootnote\relax\footnote{Keywords: Exponential-Algebraic Closedness, abelian varieties, homology.}
\let\thefootnote\relax\footnote{2020 MSC: 03C65, 14K12, 14K20.}

It is well-known that subsets of the complex numbers that are definable in the language of rings are finite or cofinite; in fact, this is true for every structure that is elementarily equivalent to $\mathbb{C}$. This property is known as \textit{strong minimality}, and the structures which satisfy it have a very tame model theory.

Let us now consider the structure $\mathbb{C}_{\exp}:=\langle \mathbb{C}, +, \cdot, 0,1, \exp \rangle$ obtained by adding to the language of rings a symbol for the complex exponential function. Strong minimality then immediately fails: the kernel of $\exp$ becomes a (countably infinite) definable set. We may however ask to what extent we move away from strong minimality: for example, if the real numbers were definable in this language then $\mathbb{C}_{\exp}$ would interpret full second-order arithmetic, and its model theory would be wild.

These questions were first looked at and formalized by Zilber, who introduced the notion of a \textit{quasiminimal structure}.

\begin{defn}
	Let $M$ be a model-theoretic structure in a countable language. We say $M$ is \textit{quasiminimal} if all definable subsets of $M$ are countable or cocountable.
\end{defn}

It is then natural to ask whether $\mathbb{C}_{\exp}$ is a quasiminimal structure. This was the content of Zilber's \textit{Quasiminimality Conjecture} (see \cite{Zil97}, \cite{Zil05}). 

Zilber's strategy to attack this problem was to define a (non-elementary) class of quasiminimal structures which is uncountably categorical (it contains exactly one structure for each uncountable cardinal up to isomorphism). The problem then became to determine whether $\mathbb{C}_{\exp}$ is isomorphic to the unique model in this class of cardinality $2^{\aleph_0}$, usually denoted by $\mathbb{B}$. 

Bays and Kirby showed in \cite{BK18} that proving that one of the axioms of Zilber's class of structures, \textit{Exponential-Algebraic Closedness}, holds in the complex numbers is actually sufficient to establish the Quasiminimality Conjecture. This axiom predicts sharp sufficient conditions for systems of equations involving polynomials and exponentials to be solvable; it is essentially a problem in complex geometry, asking to show that certain complex algebraic varieties intersect the graphs of the exponential functions: therefore, motivated by the Quasiminimality Conjecture, over the last few years several results concerning this problem have been established, see \cite{AKM}, \cite{BM}, \cite{DFT21}, \cite{Gal22} \cite{K19}, \cite{MZ}. 

Zilber's work on the model theory of $\exp$ was generalized by Bays and Kirby in \cite{BK18} to the exponential maps of abelian varieties. Given an abelian variety $A$ of dimension $g$, there is a universal covering map $\exp_A: \mathbb{C}^g \twoheadrightarrow A$. This shares many properties with the complex exponential function (in particular, for example, its kernel is countable), and hence it makes sense to consider analogues of the Quasiminimality Conjecture in this setting.

For example, given an elliptic curve $E$, it has an exponential function of the form $\exp_E:\mathbb{C} \twoheadrightarrow E$ mapping $z$ to the point $\sqvecttwo{1}{\wp(z)}{\wp'(z)}$ for some Weierstrass $\wp$-function $\wp$ and its derivative $\wp'$: is the structure $\langle \mathbb{C}, +, \cdot, 0, 1, \wp, \wp' \rangle$ quasiminimal?

As in the case of the complex exponential function, Bays and Kirby reduce the problem of quasiminimality of the structure to a question on solutions of systems of equations. In the setting of abelian varieties, Exponential-Algebraic Closedness takes the following form.

\begin{conj}[Abelian Exponential-Algebraic Closedness, see {\autocite[Sections 8 and 9]{BK18}}]
	Let $A$ be a complex abelian variety of dimension $g$, $\exp_A:\mathbb{C}^g \twoheadrightarrow A$ its exponential map, $\Gamma_{\exp_A}$ the graph of $\exp_A$.
	
	If $V$ is a free and rotund algebraic subvariety of $\mathbb{C}^g \times A$, then $V \cap \Gamma_{\exp_A} \neq \varnothing$.
\end{conj}

The notions of \textit{freeness} and \textit{rotundity} of an algebraic variety, the technical conditions of the conjecture, will be defined in Section \ref{abeac}; they are supposed to avoid any cases in which the variety ``obviously'' does not intersect $\Gamma_{\exp_A}$.

In this paper we make progress on this line of research, showing a case of Abelian Exponential-Algebraic Closedness. Specifically, we prove the following.

\begin{thm*}[Theorem \ref{abmain}]
Let $A$ be a complex abelian variety of dimension $g$ with exponential map $\exp_A:\mathbb{C}^g \rightarrow A$. Let $L \leq \mathbb{C}^g$ be a linear subspace and $W \subseteq A$ an algebraic variety such that the variety $L \times W$ is free and rotund.

Then $L \times W \cap \Gamma_{\exp_A} \neq \varnothing$.
\end{thm*}

Here $\Gamma_{\exp_A}$ denotes again the graph $$\{(z,a) \in \mathbb{C}^g \times A \mid a=\exp_A(z) \}$$ of the exponential map of $A$. 

Theorem \ref{abmain} improves a previous result of the author, \autocite[Theorem~2.12]{Gal}, which imposed stronger conditions on the space $L$.

We remark that results similar to Theorem \ref{abmain} in the context of the complex exponential function were proved by Zilber in \cite{Zil02} and in \cite{Zil11}, and by the author of the present paper in the preprint \cite{Gal22}. While the formulation of the theorems is similar, however, the proofs are quite different: in the setting of the complex exponential we use tropical geometry and the theory of amoebas, which are not available in the case of abelian varieties as these are not toric. Instead, here we use arguments from homology theory. 

The paper is structured as follows. In Section \ref{geop} we review some background material on homology and o-minimality. In Section \ref{abeac} we discuss the origin and the consequences of the Abelian Exponential-Algebraic Closedness Conjecture. In Section \ref{abmainsection} we prove the main theorem of the paper, Theorem \ref{abmain}.

\textbf{Acknowledgements.} The author wishes to thank his supervisor, Vincenzo Mantova, for suggesting to work on this problem and for many valuable comments, and the examiners of his PhD thesis, Andrew Brooke-Taylor and Jonathan Pila, for their helpful remarks on this work and its exposition. The author also wishes to thank the referees who helped catching some mistakes and simplifying the paper. This research was done as part of the author's PhD project, supported by a scholarship of the School of Mathematics at the University of Leeds.

\section{Geometric Preliminaries}\label{geop}

\subsection{Notation and Conventions}
In this paper, a complex abelian varieties $A$ of dimension $g$ will be seen as a complex algebraic variety embedded into some projective space $\mathbb{P}^N(\mathbb{C})$. The Lie algebra of $A$ will be identified with the linear space $\mathbb{C}^g$. Hence, the tangent bundle of $A$ (and its subvarieties) will be seen as the variety $\mathbb{C}^g \times A$ embedded into $\mathbb{C}^g \times \mathbb{P}^N(\mathbb{C})$.

Moreover, for any abelian subvariety $B$ of $A$ we will identify the Lie algebra $LB$ of $B$ with a subspace of $LA \cong \mathbb{C}^g$, and the Lie algebra $L(A/B)$ of the quotient abelian variety with a linear space $\mathbb{C}^{g-\dim B}$.

We will identify all these varieties with the sets of their complex points, and assume them all to be irreducible even when not explicitly stated.

\subsection{Homology and Cohomology}\label{homcohom}

The homology theory we are going to use is singular homology. We refer the reader to \autocite[Chapter~2]{Hat} for a general introduction to the subject, and to \autocite[Chapter~1]{BL} for the case of abelian varieties.

In the case of abelian varieties (and in general in the case of complex tori), the structure of the homology groups is easy to compute: in an abelian variety $A$ of dimension $g$, the $n$-th homology group is a free abelian group of rank $\vecttwo{2g}{n}$. 

As is the case in general with projective algebraic varieties, it is possible to triangulate complex algebraic subvarieties of an abelian varieties (see for example \autocite[Section~2]{Hiro} for a proof). We will therefore treat algebraic subvarieties of complex abelian varieties as cycles: given a variety $W$ and a triangulation $\mathcal{T}$, we can identify the variety with the cycle $\sum_{\sigma \in \mathcal{T}} \sigma.$ The choice of triangulation is only unique up to homological equivalence, but this is sufficient for our purposes.

We are interested in integration of certain differential forms on chains; in particular, we will study decomposable differential forms (i.e$.$ forms of degree $d$ obtained as the exterior product of $d$ forms of degree 1) attached to linear spaces. We will be interested in when the integral of these forms is non-zero, rather than in the forms themselves, and therefore it is sufficient to consider them up to scalar multiplication.

\begin{defn}
    Let $T \leq \mathbb{R}^n$ be a hyperplane defined by the equation $r_1x_1+\dots+r_nx_n=0$. We denote by $\omega_T$ (any scalar multiple of) the differential form $r_1dx_1+\dots+r_ndx_n$.
\end{defn}

We want to extend this definition to spaces of arbitrary codimension, but in that case we have more freedom in the choice of equations. The following lemma shows that this is not an issue.

\begin{lem}\label{welldef}
    Let $T_1,\dots,T_d,T_1',\dots,T_d' \leq \mathbb{R}^n$ be hyperplanes, and assume that $\bigcap_{i=1}^d T_i=\bigcap_{i=1}^{d'} T_i'$. Then there is $a \in \mathbb{R}$, $a \neq 0$, such that $$\bigwedge_{i=1}^d \omega_{T_i'}=a\bigwedge_{i=1}^d \omega_{T_i}.$$
\end{lem}

\begin{proof}
    For any $i=1,\dots,d$, the equation defining $T_i'$ is an $\mathbb{R}$-linear combination of the equations defining $T_1,\dots,T_d$; as a consequence, for all $i$ there exist $r_{i1},\dots,r_{id}$ such that $\omega_{T_i'}=\sum_{j=1}^d r_{ij}\omega_{T_j}$. Let $R$ denote the square matrix $(r_{ij})_{1\leq i,j \leq d}$: this is the matrix of a change of basis in the space of equations for $\bigcap_{i=1}^d T_i$, and therefore its determinant is non-zero. A straightforward calculation then shows that $$\bigwedge_{i=1}^d \omega_{T_i'}=\bigwedge_{i=1}^d \sum_{j=1}^d r_{ij} \omega_{T_i}=\det(R) \bigwedge_{i=1}^d \omega_{T_i}.$$
\end{proof}

\begin{defn}
    Let $T \leq \mathbb{R}^n$ be any linear subspace; let $T_1\dots,T_d$ be hyperplanes such that $\bigcap_{i=1}^d T_i=T$. We denote by $\omega_T$ (any scalar multiple of) the differential form $\bigwedge_{i=1}^d \omega_{T_i}$.
\end{defn}

Lemma \ref{welldef} implies that this is a good definition.

We will also need to consider holomorphic differential forms attached to complex subspaces of $\mathbb{C}^g$.

\begin{defn}
    Let $L \leq \mathbb{C}^g$ be a complex hyperplane, defined by the equation $l_1z_1+\dots+l_gz_g=0$. Then we denote by $\omega_L^\hol$ (any complex scalar multiple of) the differential form $l_1dz_1+\dots+l_gdz_g$.
\end{defn}

The same argument as Lemma \ref{welldef} shows that the following is a good definition.

\begin{defn}
    Let $L \leq \mathbb{C}^g$ be a complex space, and let $L_1,\dots,L_d \leq \mathbb{C}^g$ be complex hyperplanes such that $L=\bigcap_{i=1}^d L_i$. Then we denote by $\omega_L^\hol$ (any scalar multiple of) the differential form $\bigwedge_{i=1}^d \omega_{L_i}$.
\end{defn}

Thus $\omega_L^\hol$ is a holomorphic differential form of bidegree $(d,0)$.

By \autocite[Proposition 1.3.5]{BL} every complex differential form on the abelian variety $A$ is equivalent up to cohomology to a complex differential form which is invariant under translation by elements of $A$; taking the pullback of this uniquely determines a translation-invariant differential form on $\mathbb{C}^g$. We adopt the same convention, so for a complex coordinate function $v$ on $\mathbb{C}^g$, $dv$ will denote both a form on $\mathbb{C}^g$ and a form on $A$ seen as the quotient of $\mathbb{C}^g$ by a lattice.

We conclude this section with two lemmas on the relation between integrals of differential forms and transversality of intersections. Here and in the rest of the paper homological cycles in complex abelian varities will be represented by real analytic sets, so that we may discuss their real dimension.

\begin{lem}\label{nonzeroint}
Let $T \leq \mathbb{C}^g$ be a real vector subspace of dimension $2g-d$ with associated differential form $\omega_T$, and $C \subseteq A$ be a cycle of (real) dimension $d$.

If $\int_C \omega_T \neq 0$, then there is a smooth point $c \in C$ such that $T_cC \cap T=\langle 0 \rangle.$
\end{lem} 

\begin{proof}
We may fix real coordinate functions $x_1,\dots,x_{2g}$ so that $T$ is defined by the equations $x_1=\dots=x_d=0$, and hence $\omega_T$ is the form $dx_1\wedge \dots\wedge dx_d$. Thus if $\int_C \omega_T \neq 0$, there must be a smooth point $c$ of $C$ such that the linear form $(dx_1 \wedge \dots \wedge dx_d)_c:\bigwedge_{i=1}^d T_cC \rightarrow \mathbb{R}$ is non-zero. Hence, there is a basis $\{c_1,\dots,c_d\}$ of $T_cC$ such that $c_1,\dots,c_d$ are linearly independent over $T$. Hence, $T_cC+T=\mathbb{C}^g$, and as their dimensions are complementary $T_cC \cap T = \langle 0 \rangle$.
\end{proof}

Proposition \ref{nonzeroint} has a partial converse, in the case in which the cycle $C$ is a complex subvariety of $A$ and the linear space is complex.

\begin{prop}\label{nonzerointconv}
Let $L \leq \mathbb{C}^g$ be a complex vector subspace of (real) dimension $2g-2d$ with associated holomorphic differential form $\omega_L^\hol$, and $W \subseteq A$ a complex algebraic subvariety of (complex) dimension $d$.

If there is one smooth point $w \in W$ with $L \cap T_wW=\langle 0 \rangle$, then $\int_W \omega_L^\hol \wedge \o{\omega_L^\hol} \neq 0$.
\end{prop}

\begin{proof}
	Since $W$ is a complex variety, at each of its smooth points the tangent space $T_wW$ is a complex vector space of complex dimension $d$. As we noted, the form $\omega_L^\hol$ is a decomposable $(d,0)$-form, and hence $\omega_L^\hol \wedge \o{\omega_L^\hol}$ is a scalar multiple of a strongly positive form. 
	
    If there is one point $w \in W$ with $L \cap T_wW=\langle 0 \rangle$, then the pullback $\iota^*(\omega_L^\hol)$ of $\omega_L^\hol$ to $W$ is a non-zero holomorphic $(d,0)$-form on a complex projective variety of dimension $d$. Hence, $$\int_W \omega_L^\hol \wedge \o{\omega_L^\hol}=\int_W \iota^*(\omega_L^\hol) \wedge \o{\iota^*(\omega_L^\hol)} \neq 0.$$ 
\end{proof}

\subsection{Definability and O-Minimality}

We conclude this section by introducing some background on o-minimal geometry. We are going to use these facts in the proof of the main result, which needs definability of Hausdorff limits. As we expect the average reader of this article to have a background in model theory and thus be probably already familiar with o-minimality, we will not get too much into detail. The reader who wants to know more about o-minimal structures, not just on the real numbers, should consult the book \cite{vdD}. 

We follow Pila's approach from \cite{Pil} to introduce o-minimal structures on the reals.

\begin{defn}
An \textit{o-minimal structure on $\mathbb{R}$ expanding the real field} is a collection $\{\mathcal{S}_n \mid n \in \omega \}$, where each $\mathcal{S}_n$ is a set of subsets of $\mathbb{R}^n$, which satisfies the following conditions:
\begin{itemize}
\item[1.] Each $\mathcal{S}_n$ is a boolean algebra;
\item[2.] $\mathcal{S}_n$ contains every semi-algebraic subset of $\mathbb{R}^n$;
\item[3.] If $A \in \mathcal{S}_n$ and $B \in \mathcal{S}_m$, then $A \times B \in \mathcal{S}_{n+m}$;
\item[4.] If $m \geq n$, $\pi:\mathbb{R}^m \rightarrow \mathbb{R}^n$ is the projection on the first $n$ coordinates, and $A \in \mathcal{S}_m$, then $\pi(A) \in \mathcal{S}_n$;
\item[5.] The boundary of every subset of $\mathcal{S}_1$ is finite.
\end{itemize}
\end{defn}

O-minimal structures have played an important role in pure and applied model theory over the past thirty years: for example, in the aforementioned paper \cite{Pil} they were used to prove the Andr\'e-Oort Conjecture. The structure which we will use is the structure $\mathbb{R}_{\textnormal{an}}$.

\begin{defn}
	The structure $\mathbb{R}_{\textnormal{an}}$ is the smallest o-minimal structure on $\mathbb{R}$ expanding the real field which contains every globally subanalytic subset of $\mathbb{R}^n$ for each $n$.
\end{defn}

As this is the only structure we are going to refer to, we will say that a set is \textit{definable}, rather than $\mathbb{R}_{\textnormal{an}}$-definable, if it lies in the structure.

The fact that this structure exists is due to Denef and van den Dries (\autocite[Section~4]{DvdD}). It is obviously a structure which is suitable to talk about complex abelian varieties, as $\exp_A$ is an analytic map on each neighbourhood of a fundamental domain and hence definable there. Note that then, in particular, we can treat closed subgroups of abelian varieties $A$ as definable sets, as they are images under $\exp_A$ of the intersection of a linear space with finitely many fundamental domains. 

Another important feature of definable set is that they admit triangulations. The following is known as the \textit{Triangulation Theorem} for definable sets.

\begin{thm}[{\autocite[Theorem~8.1.7]{vdD}}]\label{deftri}
Every compact definable set admits a triangulation by definable sets.
\end{thm}

Therefore we may treat definable sets as cycles, and integrate differential forms on them. 

We are interested in Hausdorff limits of definable families of sets. By \textit{Hausdorff limit} we mean the limit in the topology defined by the Hausdorff distance, which we now define.

\begin{defn}
Let $A,B \subseteq \mathbb{R}^n$ be closed sets. The \textit{Hausdorff distance} between $A$ and $B$, denoted $\dhaus(A,B)$, is defined by $$\dhaus(A,B)=\max\{\sup_{a \in A}d(a,B), \sup_{b \in B}d(b,A) \}$$ where $d(x,Y)$ denotes the usual Euclidean distance between the point $x$ and the set $Y$.
\end{defn}

The reader may verify that the topology induced by this distance is $T_2$, and that therefore limits are unique.

We will call the \textit{Hausdorff limit} of a family of sets the limit of that family in the Hausdorff topology.

\begin{defn}
A definable family of definable sets is a family of sets of the form $$\{S(x) \mid x \in X  \}$$ where $S \subseteq \mathbb{R}^{m+n}$, $X \subseteq \mathbb{R}^n$ and $$S(x):=\{s \in \mathbb{R}^n \mid (s,x) \in S \}$$ are all definable sets.
\end{defn}

The tameness properties of o-minimal structures allow us to avoid pathological situations when taking Hausdorff limits of such definable collections. 

The following result follows from a theorem of Marker and Steinhorn on definable types in o-minimal theories (\autocite[Theorem~2.1]{MaSt}), which was strengthened by Pillay in \autocite[Corollary~2.4]{Pil94}. A version for semialgebraic sets was established by Br\"ocker in \autocite[Corollary~2.8]{Bro}, while a direct geometric proof appears in \autocite[Corollary~2]{KPV}. We refer the reader also to \autocite[Theorem, p.377]{LS} and \autocite[Theorem~3.1 and Proposition~3.2]{vdDL}.

\begin{thm}\label{deflim}
	Let $\{S(x) \mid x \in X\}$ be a definable family of definable compact sets over the set $X \subseteq \mathbb{R}^n$. Suppose $\gamma:[0,1) \rightarrow X$ is a definable function. Then $$S(1):=\lim_{t \rightarrow 1} S(\gamma(t))$$ exists, is a compact definable set, and $$\dim S(1) \leq \lim_{t \rightarrow 1} \dim S(\gamma(t)).$$
\end{thm}

We use Theorem \ref{deflim} to establish a property of intersection homology in the o-minimal setting.

\begin{lem}\label{inthom}
Let $A$ be an abelian variety, $W \subseteq A$ an algebraic subvariety and $\mathbb{T}$ a closed subgroup of $A$. Denote by $\{W\}$ and $\{\mathbb{T}\}$ the homology classes of $W$ and $\mathbb{T}$, and by $\{W\} \cdot \{\mathbb{T}\}$ their product in the sense of Goresky-Macpherson (see \cite{GM}).

Then $W \cap \mathbb{T}$ contains a cycle which lies in the homology class $\{W\} \cdot \{\mathbb{T}\}$.
\end{lem}

\begin{proof}
If $\mathbb{T}=A$ then $\{W\} \cdot \{A\}=\{W\}$, so the lemma holds. Hence we assume $\mathbb{T}$ is a proper subgroup.

If $W$ and $\mathbb{T}$ are dimensionally transverse (which in the smooth case coincides with the usual notion of transversality) then the lemma holds by the Goresky-Macpherson definition of the intersection pairing (see \autocite[Theorem 1]{GM}). This is the generic situation: taking a stratification of $W$ into smooth sets, we may find an open dense subset $O$ of $A$ such that for every $a \in O$ the intersection $a+\mathbb{T} \cap W$ is dimensionally transverse.
	
Let $\gamma:[0,1] \rightarrow A$ be a definable continuous function with $\gamma([0,1)) \subseteq O$ and $\gamma(1)=0_A$. We can choose $\gamma$ so that $\gamma(c_1)+\mathbb{T} \neq \gamma(c_2)+\mathbb{T}$ for all $c_1 \neq c_2 \in O$: it suffices to take $\gamma:[0,1] \rightarrow O+\mathbb{T} \subseteq A/\mathbb{T}$, which is still a dense open set, and then lift it to a curve in the original space.
	
For any $c \in (0,1)$, let $$S_c:=\bigcup_{0 \leq t \leq c} (\gamma(t)+\mathbb{T}) \cap W.$$ Since all intersections $(\gamma(t)+\mathbb{T}) \cap W$ are dimensionally transverse and pairwise disjoint, it is clear that $$\dim_{\mathbb{R}} S_c= \dim_{\mathbb{R}} \mathbb{T}+\dim_{\mathbb{R}} W -2g+1$$ for every $c \in (0,1)$, and that denoting by $B_c$ the intersection $\gamma(c)+\mathbb{T} \cap W$ for $c \in [0,1)$, $\partial S_c=B_0-B_c$ as a chain.
	
Let now $S_1=\lim_{c \rightarrow 1} S_c$ and $B_1=\lim_{c \rightarrow 1} B_c$. These are both Hausdorff limits of definable families, and therefore by Theorem \ref{deflim} the dimensions do not increase: thus $\dim_\mathbb{R} S_1 \leq \dim_{\mathbb{R}} S_c$, and as $S_1$ contains every $S_c$ the dimensions need to be equal; $\dim_{\mathbb{R}} B_1 \leq \dim_{\mathbb{R}} B_c$, and a priori it could be that $\dim_{\mathbb{R}} B_1 < \dim_{\mathbb{R}} B_c$ (although a consequence of this lemma is that this is only possible if the $B_c$'s are homologically trivial). Both $S_1$ and $B_1$ admit a triangulation by Theorem \ref{deftri}.
	
\textbf{Claim}: $\partial S_1=B_0-B_1$.
	
\textbf{Proof of Claim}: Consider the definable set $$\{(s,c) \mid s \in S_c \}$$ whose fibre above each $c$ is $S_c$. By a consequence of the Trivialization Theorem for o-minimal structures (see 9.2.1 of \autocite[{}~9.2.1]{vdD}) there is a simplicial complex $\mathcal{K}$ such that (after reparametrizing $\gamma$ to change the starting point if necessary) for all $c \in (0,1)$, $S_c$ is definably homeomorphic to $\mathcal{K}$. Consider the resulting simplices $\sigma_{j,c}:\Delta_{j} \rightarrow S_c$ (so each $\Delta_{j}$ is a face of the standard simplex and each $\sigma_{j,c}$ is a continuous function). For each $j$, we can define $\sigma_{j,1}:\Delta_{j} \rightarrow S_1$ by $\sigma_{j,1}(x)=\lim_{c \rightarrow 1}\sigma_{j,c}(x)$.

The boundary map on simplices preserves limits, in the sense that $\lim_{c \rightarrow 1}\partial \sigma_{j,c}=\partial \sigma_{j,1}$: this is clear as $\Delta_{j}$ is compact and $\sigma_{j,c}$ is continuous. Therefore, $$\partial S_1=\sum_{j=1}^{M} \partial \sigma_{j,1}=\lim_{c \rightarrow 1} \sum_{j=1}^{M} \partial \sigma_{j,c}=\lim_{c \rightarrow 1} \left(B_0-B_c\right)=B_0-B_1$$ so we are done. 

This proves the claim.

Then $B_1$ is homologous to $B_0$, which is a dimensionally transverse intersection and therefore lies in the correct homology class.
\end{proof}

\begin{rem}
Lemma \ref{inthom} fails without the o-minimality assumption: after \autocite[Theorem~VI.11.10]{Bre} an example is discussed in which the intersection of cycles $A$ and $B$ does not contain any cycle which lies in the product homology class, although every neighbourhood of the intersection does. In that case a curve of $\sin \left( \frac 1x\right)$-type is used: this is a highly non-definable object.
\end{rem}

\section{Abelian Exponential-Algebraic Closedness}\label{abeac}

\subsection{The Conjecture}

In this subsection we give a brief account of the origin of the Exponential-Algebraic Closedness Conjecture for abelian varieties.

The question is inspired by an analogous problem for the complex exponential function, posed by Zilber in \cite{Zil05}. In turn, this originated from Zilber's work on the model theory of the exponential function, and in particular from the \textit{Quasiminimality Conjecture}.

\begin{conj}[Quasiminimality Conjecture, \cite{Zil97}]\label{quaconj}
Let $\mathbb{C}_\exp$ denote the structure on the complex numbers in the language $\mathcal{L}=\{+,-,\cdot,1,0,\exp \}$ of fields equipped with an exponential. 

Let $A \subseteq \mathbb{C}$ be a finite set. Every subset of $\mathbb{C}$ which is invariant under all automorphisms of $\mathbb{C}_\exp$ which fix $A$ is countable or cocountable (we say that the structure $\mathbb{C}_\exp$ is \textit{quasiminimal}).
\end{conj}

Conjecture \ref{quaconj} stemmed out of a belief that all the definable subsets of $\mathbb{C}$ in the structure $\mathbb{C}_\exp$ that are not finite or cofinite should in some sense originate from the fact that $\exp$ has countable fibres. If the Quasiminimality Conjecture were true, then definable subsets of powers of $\mathbb{C}$ would have good tameness properties, in particular good geometric properties similar to those of complex algebraic varieties.

Further work of Zilber and of Bays-Kirby (see \cite{Zil05}, \cite{BK18}) led to reducing the Quasiminimality Conjecture to a consequence of another conjecture on the complex numbers, called \textit{Exponential-Algebraic Closedness}. We do not state this conjecture here, but we simply say that it predicts sufficient conditions for an algebraic subvariety of $\mathbb{C}^n \times (\ctimes)^n$ to intersect the graph of (the $n$-th Cartesian power of) the exponential function. In other words, it is a conjecture on solvability of systems of equations which involve polynomials and the exponential map.

In fact, the work of Bays and Kirby applies to more general settings than the exponential function, and it encapsulates model-theoretic properties of the exponential maps of abelian varieties. This was done using the construction of $\Gamma$-\textit{fields} (see \autocite[Section~3]{BK18}): these are fields together with a predicate for a certain module which mimics the behaviour of the graph of a transcendental group homomorphism such as $\exp$ or the exponential of an abelian variety. One of the main results of \cite{BK18} (Theorem 1.7 there) is that $\Gamma$-fields give rise to quasiminimal structures, and the goal is to show that the corresponding structures on the complex numbers are isomorphic to $\Gamma$-fields.

The reader might now guess that the Exponential-Algebraic Closedness Conjecture for abelian varieties (also known as Abelian Exponential-Algebraic Closedness Conjecture) that we are about to state plays the same role as the correspondent statement for the complex exponential: it predicts sufficient conditions for solvability of systems which involve algebraic operations and the exponential of an abelian variety, and its main consequence would be model-theoretic (quasiminimality of a structure).

With this in mind, we begin to lay the ground for the statement of the conjecture. The notions to define are then \textit{freeness} and \textit{rotundity} of an algebraic subvariety of $\mathbb{C}^g \times A$, which should work as the sufficient conditions for the variety to intersect the graph of the exponential map.

\begin{defn}\label{freerotundab}
	Let $A$ be an abelian variety of dimension $g$, $V \subseteq \mathbb{C}^g \times A$ an algebraic variety.
	
	Let $\pi_1$ and $\pi_2$ denote the projections of $\mathbb{C}^g \times A$ on $\mathbb{C}^g$ and $A$ respectively.
	
	The variety $V$ is \textit{free} if $\pi_1(V)$ is not contained in a translate of the Lie algebra of a non-trivial abelian subvariety of $A$ and $\pi_2(V)$ is not contained in a translate of a non-trivial abelian subvariety of $V$.
	
	Given an abelian subvariety $B \leq A$, denote by $\pi_B$ the projection $\pi_B:\mathbb{C}^g \times A \rightarrow \mathbb{C}^{g-\dim B} \times A/B$.
	
	The variety $V$ is \textit{rotund} if $\dim \pi_B(V) \geq \dim A/B$ for every abelian subvariety $B$.
\end{defn}

Note that since this definition quantifies over abelian subvarieties it becomes almost trivial in the case of a simple abelian variety. If $A$ is simple, then $V \subseteq \mathbb{C}^g \times A$ is free if and only if $\pi_1(V)$ and $\pi_2(V)$ are both positive-dimensional, and it is rotund if and only if $\dim V \geq g$.

The conjecture is then that every free rotund algebraic variety intersects the graph of the exponential function. While it is not stated there in this form, it is implicit in \autocite[Sections~8~and~9]{BK18}. 

\begin{conj}[Abelian Exponential-Algebraic Closedness]\label{abeacconj}
	Let $A$ be an abelian variety of dimension $g$, $\exp_A:\mathbb{C}^g \rightarrow A$ its exponential map, $\Gamma_{\exp_A}$ the graph of $\exp_A$.
	
	For every free and rotund algebraic variety $V \subseteq \mathbb{C}^g \times A$, $V \cap \Gamma_{\exp_A} \neq \varnothing$.
\end{conj}

As we mentioned, the main interest in the Abelian Exponential-Algebraic Closedness Conjecture lies in the fact that if it were confirmed it would imply quasiminimality of a structure on the complex numbers.

\begin{thm}[{\autocite[Theorem~1.9]{BK18}}]
	Let $A$ be a simple abelian variety. If Conjecture \ref{abeacconj} holds for all powers of $A$, then the structure $(\mathbb{C}, +, - \cdot, 0,1, \Gamma_{\exp_A})$ is quasiminimal.
\end{thm}

Therefore, just like Exponential-Algebraic Closedness, Conjecture \ref{abeacconj} has a model-theoretic motivation which ties it to the theory of quasiminimal structures.

\subsection{Consequences of Freeness and Rotundity}

In this section we gather some geometric consequences of the properties freeness and rotundity. In particular, we are going to work on the following complex analytic function.

\begin{defn}
Let $V \subseteq \mathbb{C}^g \times A$ be an algebraic variety. The \textit{$\delta$-map} of the variety $V$ is the function $\delta:V \rightarrow A$ defined by $(v_1,v_2) \mapsto v_2-\exp_A(v_1)$.
\end{defn}

The following Lemma establishes a connection to between freeness and rotundity and openness of the $\delta$-map of a variety. It is similar to a result of Kirby for the exponential function (\cite[Proposition~6.2~and~Remark~6.3]{K19}, the proof idea is essentially the same), and a less general version was present in the author's earlier work, see \cite[Proposition~2.11~and~Theorem2.12]{Gal}.

\begin{lem}\label{abopen}
	Let $V \subseteq \mathbb{C}^g \times A$ be a rotund variety, and let $\delta:V \rightarrow A$ map $v=(v_1,v_2)$ to $v_2 - \exp_A(v_1)$. 
	
	Then there is a Zariski-open dense subset $V^\circ \subseteq V$ such that for every $v \in V$, $\delta$ is open at $v$. Moreover, if $A$ is simple, then $V^\circ=V$.
\end{lem}

\begin{proof}
	Let $V \subseteq \mathbb{C}^g \times A$ and $\delta$ be as in the statement. $\delta$ is clearly a complex analytic map, and therefore its fibres are complex analytic sets: if there is $a \in A$ such that $\delta^{-1}(a)$ is a complex analytic subset of $V$ of dimension $\dim V -g$, then the image of $\delta$ contains, around $a$, a $g$-dimensional complex analytic set, and in particular the map is open on some neighbourhood of any point $v \in \delta^{-1}(a)$.
	
	The fibre $\delta^{-1}(a)$ coincides with the set $$\{v \in V \mid \delta(v)=a \}=\{(v_1,v_2) \in V \mid v_2-\exp_A(v_1)=a \}=$$ $$=\{(v_1,v_2) \in V \mid \exp(v_1)=v_2-a \}$$ and therefore it is a translate of the intersection $\Gamma_{\exp_A} \cap (V-(0,a)).$
	
	Consider the family of varieties $\mathcal{V}:=\{V-(0,a) \mid a \in A\}$. By Kirby's Uniform Ax-Schanuel Theorem \cite[Theorem~4.3]{K09} there is a finite collection of abelian subvarieties $\mathcal{C}$ such that any irreducible component of any intersection $(V-(0,a)) \cap \Gamma_{\exp_A}$ of dimension larger than $\dim V - g$ is contained in a translate of the tangent bundle of some $B \in \mathcal{C}$, $\gamma +(LB \times B)$. Note that if $A$ is simple this implies that there cannot be any such irreducible component, and therefore $\delta$ is open, proving the ``moreover'' part of the statement.
	
	If $B$ is minimal with this property, then $$\dim ( (\gamma + (LB \times B)) \cap (V-(0,a))) > \dim V -g + \dim B.$$ To see this, consider the variety $V_B:=(LB \times B) \cap (-\gamma + V-(0,a))$: as a subvariety of $LB \times B$, and by minimality of $B$, it must satisfy $$ \dim (V_B \cap \Gamma_{\exp_B})=\dim V_B - \dim B$$ as if it did not we could apply the uniform Ax-Schanuel Theorem in $B$ and obtain a smaller abelian subvariety. As $\dim (V_B \cap \Gamma_{\exp_B}) > \dim V - g$, we have $\dim V_B > \dim V - g + \dim B$.  
	
	By rotundity, if $\pi_B$ denotes the quotient $\mathbb{C}^g \times A \twoheadrightarrow \mathbb{C}^{g-\dim B} \times A/B$ we must have $$\dim \pi_B(V) \geq g - \dim B$$ and therefore for almost all $\gamma \in\mathbb{C}^g \times A$, $$\dim  ((\gamma+(LB \times B)) \cap (V-(0,a))) = \dim V - g + \dim B. $$ Therefore after removing a Zariski-closed proper subset of $V$ we obtain a set $V^\circ \subseteq V$, Zariski-open dense in $V$, such that for every point in $V$ and every abelian subvariety $B$ we must have that $$\dim  ((\gamma+(LB \times B)) \cap (V-(0,a))) = \dim V - g + \dim B.$$ By the argument above, this implies that the $\delta$-map of $V$ is open at every point of $V^\circ$: if it were not, then the dimension equality would not be satisfied for some abelian subvariety $B \subseteq A$.
\end{proof}

When the variety has the form $L \times W$ we know more about the structure of the set $V^\circ$.

\begin{prop}
	If $L \leq \mathbb{C}^g$ is a linear space, $W \subseteq A$ is algebraic, and $L \times W$ is rotund, then there is a Zariski-open dense subset $W^\circ$ of $W$ such that $\delta$ is open at every point of $L \times W^\circ$.
\end{prop}

\begin{proof}
Suppose $\delta$ is open around the point $(l,w) \in L \times W$. Then there exist open neighbourhoods $U_L \subseteq L$ of $l$ and $U_W \subseteq W$ of $w$ such that the restriction of $\delta$ to $U_L \times U_W$ is open.

Now let $l' \in L$ be any other point. Since $L$ is a subgroup of $\mathbb{C}^g$, $(l'-l)+U_L$ is a neighbourhood of $l'$ in $L$. 

Let then $O$ be any open subset of $((l'-l)+U_L \times U_W)$. It is then clear that $O$ is a translate by $(l'-l,0_A)$ of an open subset $O' \subseteq U_L \times U_W$. The image of $O'$ under $\delta$ is open, because $\delta$ is open on $U_L \times U_W$, and therefore so is the image of $O$ because $\delta(O)=\delta((l-l',0_A)+O')=\delta(O')-\exp_A(l-l')$ is a translate of an open set, and hence open.

Thus if $\delta$ is open on a neighbourhood of $(l,w)$ then it is open on a neighbourhood of $(l',w)$ for any $l' \in L$. Therefore the Zariski-open dense set on which it is open must have the form $L \times W'$. 
\end{proof}

Lemma \ref{abopen} has a converse statement.

\begin{lem}\label{rotcon}
	Let $L \leq \mathbb{C}^g$ be a linear space, $W \subseteq A$ an algebraic variety and $\delta$ the $\delta$-map of $L \times W$. If there is a point $(l,w) \in L \times W$ such that $\delta$ is open at $(l,w)$, then the variety $L \times W$ is rotund.
\end{lem}

\begin{proof}
	Let $B\leq A$ be an abelian subvariety, and let $LB \cong \mathbb{C}^{g-\dim B} \leq \mathbb{C}^g$ be its tangent space at identity.
	
	Let $\pi_B:\mathbb{C}^g \times A \twoheadrightarrow \mathbb{C}^{g-\dim B} \times A/B$ denote the quotient map. Consider also the partial quotients $p:\mathbb{C}^g \twoheadrightarrow \mathbb{C}^g/LB$ and $q:A \twoheadrightarrow A/B$. We need to show that $\dim \pi_B(L \times W) \geq \dim A - \dim B$. 
	
	If $U_L \subseteq L$ and $U_W \subseteq W$ are open subsets such that $\delta$ is open on $U_L \times U_W$, then we have that $\delta(U_L \times U_W)$ is an open subset of $A$; thus $q(\delta(U_L \times U_W))$ is an open subset of $A/B$.
	
	Let $\delta_Q$ denote the $\delta$-map of the variety $p(L) \times q(W)$. 
	
	Consider $\delta_Q \circ \pi_B(l,w)$. This is $$\delta_Q(\pi_B(l,w))=\delta_Q(p(l),q(w))=q(w)-\exp_{A/B}(p(l))=$$ 
	$$=w+B-\exp_{A/B}(l+LB)=w+B-(\exp_A(l)+B)=(w-\exp_A(l))+B=$$ $$=\delta(l,w)+B=q(\delta(l,w))$$
	
	Therefore $\delta_Q \circ \pi_B=q \circ \delta$: the square 
	
	$$\begin{tikzpicture}
	\node (A) at (0,2) {$L \times W$};
	\node (B) at (4,2) {$A$};
	\node (C) at (4,0) {$A/B$};
	\node (D) at (0,0) {$p(L) \times q(W)$};
	
	\draw[->] (A) -- (B) node[midway, above]{$\delta$};
	\draw[->] (B) -- (C) node[midway, right]{$q$};
	\draw[->] (A) -- (D) node[midway, right]{$p \times q$};
	\draw[->] (D) -- (C) node[midway, above]{$\delta_Q$};
	\end{tikzpicture}$$
	
	commutes.
	
	Therefore, $\delta_Q(\pi_B(U_L \times U_W))$ is an open subset of $A/B$. This implies that $\dim (\pi_B(U_L \times U_W)) \geq \dim (A/B)$, because the image of an analytic map cannot have larger dimension than its domain. So rotundity holds.
\end{proof}

We make a useful reduction, showing that we can assume without loss of generality that that $\dim L + \dim W =g$.

\begin{lem}\label{generichypab}
	Let $L \times W$ be a free rotund variety in $\mathbb{C}^g \times A$. Then there is a space $L'' \subseteq L$ such that $L'' \times W$ is free and rotund, and $\dim L'' + \dim W =g$.
\end{lem}

\begin{proof}
We use Lemma \ref{rotcon}. Suppose that $L \times W$ is free and rotund, and that $\dim L + \dim W > g$. Since $L \times W$ is rotund, the image of the $\delta$-map of $L \times W$ has non-empty interior; let $a$ be a point in the interior of $\im(\delta)$. Then the fibre $\delta^{-1}(a)$, which is equal to the analytic set $(L \times (a+W)) \cap \Gamma_{\exp_A}$ has an irreducible component of dimension $\dim L \times W -g=\dim L + \dim W -g$ (if it did not, it could not have a $g$-dimensional neighbourhood contained in the image). Therefore, intersecting $L$ with a generic hyperplane $H \leq \mathbb{C}^g$ we can make sure that the variety $(L \cap H) \times W$ is still free, and that $((L \cap H) \times (a+W)) \cap \Gamma_{\exp_A}$ has an irreducible component of dimension $\dim (L \cap H) + \dim W -g=\dim L -1 +\dim W -g$. Therefore, $a$ is an interior point of the $\delta$-map of $(L \cap H) \times W$, which is therefore a free rotund variety of dimension $\dim L + \dim W -1$. Iterating this process finitely many times we obtain the space $L''$ such that $\dim L + \dim W=g$, such that the variety $L \times W$ is still free and rotund.
\end{proof}

Therefore we will freely assume that $\dim L = \codim W$ when we need it.

Note that we can use the characterization of rotundity in simple varieties to easily establish a partial result in that setting.

\begin{prop}\label{simplecase}
	Let $A$ be a simple abelian variety, $L \times W\subseteq \mathbb{C}^g \times A$ a free rotund variety, with $L \leq \mathbb{C}^g$ a linear space and $W \subseteq A$ a smooth algebraic variety.
	
	Then $\exp_A(L) \cap W \neq \varnothing$.
\end{prop}

\begin{proof}
	Let $\mathbb{T}=\o{\exp_A(L)}$. Let $\delta_L$ denote the $\delta$-map of $L \times W$, and $\alpha_{\mathbb{T}}: \mathbb{T} \times W \rightarrow A$ the similar map $(t,w) \mapsto w-t$. 
	
	Since $\mathbb{T}$ and $W$ are compact subsets of $A$, the image of $\alpha_{\mathbb{T}}$ is compact and hence closed. The map $\delta_L$ is open on $L \times W$ by the ``moreover'' part of Lemma \ref{abopen}.
 
	If $a \in \im(\alpha_{\mathbb{T}})$ then there are $(t,w) \in \mathbb{T} \times W$ such that $w-t=a$. Consider the variety $(t+L) \times W$: the $\delta$-map of this variety is open around $(t,w)$, and therefore its image contains a neighbourhood of $a$; this image is contained in $\im(\alpha_{\mathbb{T}})$, which therefore also contains a neighbourhood of $a$. Therefore $\im(\delta_\mathbb{T})$ is open and closed, so the map is surjective, and thus there is $w \in W \cap \mathbb{T}$. 
	
	Since $\delta_L$ is open around $(0,w)$, there are open neighbourhoods $U_0$ of $0$ in $L$ and $U_w$ of $w$ in $W$ such that $\delta_L(U_0 \times U_w)$ is an open neighbourhood of $w$ in $A$. Thus its intersection with $\mathbb{T}$ is an open neighbourhood of $w$ in $\mathbb{T}$, and it contains points of $\exp_A(L)$ (as it is dense in $\mathbb{T}$). Therefore there are points $(l',w') \in L \times W$ such that $w'-\exp_A(l') \in \exp_A(L)$, and thus $w' \in W \cap \exp_A(L)$.
\end{proof}

\section{Abelian E.A.C. for Varieties of the form $L \times W$}\label{abmainsection}

In this section we prove the main result of the paper.

\begin{thm}\label{abmain}
Let $A$ be an abelian variety of dimension $g$, $L \times W$ a free rotund variety with $L \leq \mathbb{C}^g$ a linear space and $W \subseteq A$ an algebraic variety.

Then $\exp_A(L) \cap W \neq \varnothing$.
\end{thm}

\begin{basass}\label{basic}
For the rest of the section, even when not explicitly specified we will adopt the following conventions:

\begin{itemize}
\item[(1)] $A$ is an abelian variety of dimension $g$ with exponential map $\exp_A:\mathbb{C}^g \rightarrow A$;
\item[(2)] $L \times W$ is a free rotund subvariety of $\mathbb{C}^g \times A$, with $L \leq \mathbb{C}^g$ linear, $W \subseteq A$ an algebraic subvariety, and $\codim L = \dim W=d$;
\item[(3)] $\mathbb{T}=\o{\exp_A(L)}$ is the closure of the exponential of $L$ in $A$; it is a closed subgroup of $A$, and $T$ denotes the real subspace of $\mathbb{C}^g$ such that $\exp_A(T)=\mathbb{T}$;
\item[(4)] $\delta:L \times W \rightarrow A$ is the $\delta$-map of $L \times W$, which maps $(l,w)$ to $w-\exp_A(l)$;
\item[(5)] The differential forms $\omega_L^\hol$ and $\omega_T$ are attached to the spaces $L$ and $T$ respectively according to the definitions in Subsection \ref{homcohom}
\end{itemize}
\end{basass}

We need to combine the results from Sections \ref{geop} and \ref{abeac}.

\begin{prop}
Let $L \times W$ be a free rotund variety, and assume $\dim L+\dim W=g$. Then $\int_W \omega_L^\hol \wedge \o{\omega_L^\hol} \neq 0$.
\end{prop}

\begin{proof}
By Lemma \ref{abopen}, the $\delta$-map of $L \times W$ is open on a set of the form $L \times W^\circ$ with $W^\circ \subseteq W$ Zariski-open dense. This implies that there exists at least one point $w \in W$ such that $T_wW+L=\mathbb{C}^g$. Therefore the statement holds by Proposition \ref{nonzerointconv}.
\end{proof}

We now study the interaction between the differential forms $\omega_T$ and $\omega_L^\hol$. Letting $d=\dim W = \codim L$, we know that $\omega_L^\hol$ is a holomorphic form of degree $(d,0)$.

\begin{prop}
If $T' \leq \mathbb{C}^g$ is a real vector subspace such that $L=T \cap T'$, then (up to scalar multiplication) $\omega_T \wedge \omega_{T'}=\omega_L^\hol \wedge \o{\omega_L^\hol}$.
\end{prop}

\begin{proof}
$L$ is defined by $d$ equations of the form $l_{j1}z_1+\dots+l_{jg}z_g=0$, for some $l_{j1},\dots,l_{jg} \in \mathbb{C}$, for $j=1,\dots,d$. Each of these equations can be split into its real and imaginary parts; let $R_L$ be the real subspace of $\mathbb{C}^g$ defined by the equations $\Re(l_{j1}z_1+\dots+l_{jg}z_g)=0$, for $j=1,\dots,d$, and $I_L$ the real subspace of $\mathbb{C}^g$ defined by the equations $\Im(l_{j1}z_1+\dots+l_{jg}z_g)=0$, again for $j=1,\dots,g$. Then by definition we see that $\omega_L^\hol=\omega_{R_L}+i\omega_{I_L}$. Hence, we have that $$\omega_L^\hol \wedge \o{\omega_L^\hol}=(\omega_{R_L}+i\omega_{I_L}) \wedge (\omega_{R_L} - i \omega_{I_L})=(-2i) \left( \omega_{R_L} \wedge \omega_{I_L} \right). $$ Since $L=R_L \cap I_L=T \cap T'$ we have that $\omega_{R_L} \wedge \omega_{I_L}=\omega_T \wedge \omega_{T'}$ up to multiplication by a non-zero scalar, so the statement holds.
\end{proof}

Moreover, we note that the cohomology class of $\omega_T$ and the homology class of $\mathbb{T}$ are essentially Poincar\'e duals.

\begin{prop}\label{poindual}
The cohomology class of $\omega_T$ is (a scalar multiple of) the Poincar\'e dual of the homology class of $\mathbb{T}$.
\end{prop}

\begin{proof}
This follows from \autocite[Lemma~4.10.2]{BL}.
\end{proof}

We are now ready to examine the consequences of freeness and rotundity on the homological level.

\begin{lem}\label{integral}
Let $L \times W$ satisfy the Basic Assumptions \ref{basic}, and suppose the intersection $\mathbb{T} \cap W$ is dimensionally transverse. Let $T' \leq \mathbb{C}^g$ be a real vector subspace such that $L=T \cap T'$.

Then $\int_{W \cap \mathbb{T}} \omega_{T'} \neq 0$.
\end{lem}

\begin{proof}
By the cup product-intersection duality, we know that $$\int_{W \cap \mathbb{T}} \omega_{T'}=\int_A [W] \wedge \omega_T \wedge \omega_{T'}=\int_W \omega_L^\hol \wedge \o{\omega_L^\hol} \neq 0$$ where $[W]$ is used to denote any form lying in the cohomology class that is dual to the homology class of $W$. We have used Proposition \ref{poindual} to note that $[W] \wedge \omega_T=[W \cap \mathbb{T}]$.
\end{proof}

By Lemma \ref{inthom}, $\mathbb{T} \cap W$ contains a cycle which lies in the product homology class. Therefore we have the following.

\begin{prop}\label{transcycle}
Let $L \times W$ satisfy the Basic Assumptions \ref{basic}. Then there are a cycle $C \subseteq \mathbb{T} \cap W$ and a point $c \in C$ such that $T_cC \cap T' = \langle 0 \rangle$.
\end{prop}

\begin{proof}
By Lemmas \ref{integral} and \ref{inthom}, $W \cap \mathbb{T}$ contains a cycle $C$ such that $\int_C \omega_{T'} \neq 0$. By Proposition \ref{nonzeroint}, this implies that there is $c \in C$ such that $T_cC$ and $T'$ intersect transversely.
\end{proof}

We are then ready to prove the main theorem of this paper, Theorem \ref{abmain}.

\begin{proof}[Proof of Theorem \ref{abmain}]
Since $L \times W$ is free and rotund, by Proposition \ref{transcycle} there are a cycle $C \subseteq \mathbb{T} \cap W$ and a point $c_0 \in C$ such that $T_{c_0}C + T'=\mathbb{C}^g$. Note that $$\dim_{\mathbb{R}} C=\dim_{\mathbb{R}} T+\dim_{\mathbb{R}} W -2g$$ and $$\dim_{\mathbb{R}} T'=2g-\dim_{\mathbb{R}}T+\dim_{\mathbb{R}} L$$ therefore $$\dim_{\mathbb{R}} T_{c_0} C + \dim_{\mathbb{R}} T'=\dim_{\mathbb{R}} W + \dim_{\mathbb{R}} L=2g.$$ Thus, $T_{c_0} C \cap T'=\langle 0 \rangle$. Since $L \leq T'$, this implies in particular that $T_{c_0} C \cap L=\langle 0\rangle$.

Therefore, $$\dim_{\mathbb{R}}(T_{c_0}C+L)=\dim_{\mathbb{R}} T_{c_0}C+\dim_{\mathbb{R}} L=$$ $$=\dim_{\mathbb{R}} T + \dim_{\mathbb{R}} W -2g +\dim_{\mathbb{R}} L=\dim_{\mathbb{R}} T.$$ As both spaces are contained in $T$, then, we have that $T_{c_0}C+L=T$.

This means that there is a neighbourhood $U\subseteq C$ such that $U+\exp_A(L)$ is open in $\mathbb{T}$. Therefore, since $\exp_A(L)$ is dense in $T$, there is $(l,c) \in L \times C$ such that $c-\exp_A(l) \in \exp_A(L)$; and therefore $c \in \exp_A(L) \cap C \subseteq \exp_A(L) \cap W$.
\end{proof}

Finally, we improve the result to show that actually the intersection between $\exp_A(L)$ and $W$ is Zariski-dense in $W$.

To do this, we are going to need a few preliminary results.

\begin{prop}\label{transint}
Let $L \times W$ satisfy the Basic Assumptions \ref{basic}. Then there is $w \in \exp_A(L) \cap W$ such that $T_wW+L=\mathbb{C}^g$.
\end{prop}

\begin{proof}
Theorem \ref{abmain} implies that there is $w_0 \in \exp_A(L) \cap W$ such that $\delta$ is open around $(0,w_0)$. Let $U \subseteq W$ be a neighbourhood of $w_0$ such that $\delta$ is open on $L \times U$. In every neighbourhood of $w_0$ contained in $U$ there are points of $\exp_A(L)$, by density of $\exp_A(L)$ in its closure: therefore, $\exp_A(L) \cap U$ is a countable set with no isolated points.

Assume that for every point $w \in \exp_A(L) \cap U$, $L+T_wW \neq \mathbb{C}^g$. Clearly the set of points $w \in W$ for which $L+T_wW \neq \mathbb{C}^g$ is closed: thus, since $\exp_A(L) \cap U$ is countably infinite, $w_0$ lies in a positive-dimensional set of points with this property, and therefore $\delta$ cannot have discrete fibres and be open around $(0,w_0)$. Therefore there must be a point $w$ arbitrarily close to $w_0$ with $T_wW+L=\mathbb{C}^g$: in other words, if $T_{w_0} W+L\neq \mathbb{C}^g$ but $\delta$ is open at $(0,w_0)$, then $w_0$ is isolated in the set of points with this property.
\end{proof}

Recall that a subspace $S \leq \mathbb{C}^n$ is \textit{totally real} if $S \cap iS=\langle 0 \rangle$, and that a submanifold $N$ of a complex manifold $M$ is totally real if $T_nN$ is totally real in $T_nM$ for every $n \in N$. 

We then consider the following result on complex manifolds.

\begin{prop}\label{totreal}
Let $M$ be a complex manifold with $\dim M=n$, $N \subseteq M$ a totally real submanifold with $\dim_{\mathbb{R}} N=n$, and $f:M \rightarrow \mathbb{C}$ a holomorphic function. 

If $f$ vanishes on $N$ then $f \equiv 0$ on $M$.
\end{prop}

\begin{proof}
Suppose $f$ vanishes on $N$. Consider the complex submanifold $M' \subseteq M$ defined by $\{z \in M \mid f(z)=0\}$

Then for every $n \in N$, $T_nN \subseteq T_nM'$. Since $T_nN$ is a totally real subspace of real dimension $n$, it cannot be contained in a proper complex subspace of $T_nM$: therefore, $T_nM'=T_nM$ and as a consequence $M'=M$.
\end{proof}

\begin{lem}\label{L1}
Let $L \times W$ satisfy the Basic Assumptions \ref{basic}, and let $L_1:=T+iT$ be the smallest complex subspace of $\mathbb{C}^g$ which contains $T$.

Then every holomorphic function $f:W \rightarrow A$ which vanishes on $\exp_A(L) \cap W$ vanishes on $\exp_A(L_1) \cap W$.
\end{lem}

\begin{proof}
Let $c:=\dim_{\mathbb{R}}(T \cap iT) - \dim_{\mathbb{R}} L$. Then $$\dim_{\mathbb{R}} T= \frac{\dim_{\mathbb{R}} L_1+c+\dim_{\mathbb{R}} L}{2}$$ and therefore at a point $w \in W \cap \exp_A(L)$ where $T_wW+L=\mathbb{C}^g$ (as given by Proposition \ref{transint}) we have $$\dim_{\mathbb{R}} T_wW \cap T=\dim_{\mathbb{R}}W + \dim_{\mathbb{C}} L_1+ \frac{c+\dim_{\mathbb{R}} L}{2} - 2g = $$ $$=\frac{\dim_{\mathbb{R}} W}{2}+\frac c2+\dim_{\mathbb{C}} L_1-g.$$

Note that $$\dim_\mathbb{R} (T_wW \cap T \cap iT)=\dim_{\mathbb{R}} T_wW+ \dim_{\mathbb{R}} (T\cap iT)-2g=$$ $$=\dim_{\mathbb{R}}T_wW + c + \dim_{\mathbb{R}} L -2g=c$$ and $$\dim_{\mathbb{R}} ((T_wW \cap T) + i(T_wW \cap T))=\dim_{\mathbb{R}} W + c + \dim_{\mathbb{R}} L_1  -2g-c=$$ $$=\dim_{\mathbb{R}} W + \dim_{\mathbb{R}} L_1 -2g=\dim_{\mathbb{R}} (T_wW \cap L_1)$$ so $T_wW \cap T$ is a totally real subspace of $T_wW \cap L_1$ of half its dimension, so it is not contained in any proper complex subspace. Therefore there is a neighbourhood $U \subseteq L_1$ of $0$ such that $\exp_A(U)$ is complex manifold, and $W \cap \mathbb{T}$ is at least locally a totally real submanifold of $W \cap \exp_A(U)$, of dimension $\frac{\dim_{\mathbb{R}} (W \cap U)}{2}$ (``locally'' in the sense that there is a neighbourhood $U' \subseteq W$ of $w$ such that $U' \cap \mathbb{T}$ is a totally real submanifold of the complex manifold $U' \cap \exp_A(U)$). Thus by Proposition \ref{totreal}, there is no holomorphic function which vanishes on $W \cap \mathbb{T}$ unless it vanishes on $W \cap U$, and thus on $W \cap \exp_A(L_1)$.
\end{proof}

Lemma \ref{L1} is sufficient to establish Zariski-density of $\exp_A(L) \cap W$ in $W$ in the case in which $T+iT=\mathbb{C}^g$. To prove the more general case we need an inductive argument.

\begin{thm}
Let $L \times W$ satisfy the Basic Assumptions \ref{basic}. Then $\exp_A(L) \cap W$ is Zariski-dense in $W$.
\end{thm}

\begin{proof}
If $T+iT=L_1 \lneq \mathbb{C}^g$, then let $k=k(L)$ be the length of the following chain of inclusions: $$L=L_0 \leq T=T_0 \leq L_1 \leq T_1 \leq \dots \leq L_{k-1} \leq T_{k-1} \leq L_k=\mathbb{C}^g$$ where each $L_{j+1}=T_{j}+iT_{j}$ is the smallest complex subspace of $\mathbb{C}^g$ that contains $T_j$ and each $T_{j+1}$ is the real subspace of $\mathbb{C}^g$ such that $\o{\exp_A(L_j)}=\exp_A(T_{j+1})$. Note that all the inclusions are proper, except possibly the last one: if there is a $j$ such that $T_{j-1}=L_j$ then necessarily $T_{j-1}=\mathbb{C}^g$ or the variety $L \times W$ would not be free.

Lemma \ref{L1} then shows that every holomorphic function which vanishes on $W \cap \exp_A(L)$ vanishes on $W \cap \exp_A(L_1)$: this is the $k=1$ case. Since $k$ is finite, repeated applications of the lemma show that no algebraic function can vanish on $\exp_A(L) \cap W$, which is therefore Zariski-dense in $W$.
\end{proof}

\subsection{An Example}

We conclude with a concrete example of the kind of problem that this argument applies to, with the abelian variety $A$ a product of elliptic curves. 

Recall (for example from \autocite[Section~2.3.4]{Hid}) that for every elliptic curve $E \cong \mathbb{C}/\Lambda$ there is a meromorphic function $\wp:\mathbb{C} \smallsetminus \Lambda \rightarrow \mathbb{C}$, $\Lambda$-periodic with double poles at each point of $\Lambda$ such that $E$ may be embedded in the projective space $\mathbb{P}^2(\mathbb{C})$ as the set of all points of the form $\sqvecttwo{1}{\wp(z)}{\wp'(z)}$ for $z \in \mathbb{C} \setminus \Lambda$ together with the point at infinity $\sqvecttwo{0}{0}{1}$. In other terms, the exponential map of $E$ has the form $$z \mapsto \sqvecttwo{1}{\wp(z)}{\wp'(z)}$$ (the fact that it is a group homomorphism is \autocite[Theorem~2.5.1]{Hid}). The function $\wp$ is called a \textit{Weierstrass $\wp$-function} for $E$.

Consider then the lattices $\Lambda_1=\mathbb{Z}+i\sqrt{2}\mathbb{Z}$ and $\Lambda_2=\mathbb{Z}+i\sqrt{5} \mathbb{Z}$. As $\sqrt{2}$ does not lie in the $\mathbb{Q}$-linear span of $\sqrt{5}$, the elliptic curves $E_1 \cong \mathbb{C}/\Lambda_1$ and $E_2 \cong \mathbb{C}/\Lambda_2$ are not isogenous. Let $\wp_1$ and $\wp_2$ denote the relative Weierstrass $\wp$-functions, so that for $j=1,2$ a point in $E_j$ has the form $\sqvecttwo{1}{\wp_j(z)}{\wp'_j(z)}$.

We can embed the product of the two elliptic curves in $\mathbb{P}^8(\mathbb{C})$ using the classical Segre embedding for product of projective spaces (see for example \autocite[Example~A.1.2.6(b)]{HS}) so that a point in $E_1 \times E_2 \subseteq \mathbb{P}^8(\mathbb{C})$ has the form $$\sqvecteight{1}{\wp_2(z_2)}{\wp'_2(z_2)}{\wp_1(z_1)}{\wp_1(z_1)\wp_2(z_2)}{\wp_1(z_1)\wp'_2(z_2)}{\wp'_1(z_1)}{\wp'_1(z_1)\wp_2(z_2)}{\wp'_1(z_1)\wp_2(z_2)}.$$

Consider a smooth curve $W$ contained in $E_1 \times E_2$ that is the intersection of $E_1 \times E_2$ with a hypersurface cut out in $\mathbb{P}^8(\mathbb{C})$ by a polynomial $F \in \mathbb{C}[Z_0,\dots,Z_8]$. Let moreover $$L:=\{(z_1,z_2) \in \mathbb{C}^2 \mid z_1=z_2 \}.$$ Suppose we want to find a point in $L \times W \cap \Gamma_{E_1 \times E_2}$: it corresponds to finding $z \in \mathbb{C}$ such that the corresponding point in $E_1 \times E_2$ (which we may write using the $\wp$-functions) satisfies the polynomial equation $F=0$.

To do so, let us understand what $\o{\exp_{E_1 \times E_2}(L)}$ looks like. Since, as we already noted, $\sqrt{2}$ does not lie in the $\mathbb{Q}$-linear span of $\sqrt{5}$, we have that for every fixed $z \in \mathbb{\mathbb{C}}$ the set $z+i\sqrt{2}\mathbb{Z}$ has closed image under the exponential of $E_1$, while its image under the exponential of $E_2$ looks like the set $z+i\sqrt{2}\mathbb{Z}+i\sqrt{5}\mathbb{Z}$, which is dense in $\Re(z)+i\mathbb{R}$. Therefore the closure of $\exp_{E_1 \times E_2}(L)$ is the set $$\exp(\{(z_1,z_2) \in \mathbb{C}^2 \mid \Re(z_1)=\Re(z_2) \}).$$

\begin{figure}
	\begin{center}

		\begin{tikzpicture}[scale=1.5]
		
		\draw (-2,0) -- (-1,0);
		\draw (-2,0) -- (-2,1.4);
		\draw (-1,0) -- (-1,1.4);
		\draw (-2,1.4) -- (-1,1.4);
		
		\node at (-2.15,0) {$0$};
		\node at (-2.25,1.4) {$i\sqrt{2}$};
		\node at (-1.5,0.3) {$x+iy$};
		
		\draw (1,0) -- (2,0);
		\draw (2,0) -- (2,2.2);
		\draw (1,0) -- (1,2.2);
		\draw (1,2.2) -- (2,2.2);
		
		\node at (0.85,0) {$0$};
		\node at (0.75,2.2) {$i\sqrt{5}$};
		\node at (1.5,-0.15) {$x$};
		
		\node at (-1.5, 0.5) {$\cdot$};
		\foreach \x in {0,...,11} {\node at (1.5, \x*2.2/11 ) {$\cdot$};};
		\foreach \x in {0,...,13} {\node at (1.5, \x*2.2/13 ) {$\cdot$};};
		\foreach \x in {0,...,17} {\node at (1.5, \x*2.2/17 ) {$\cdot$};};
		\foreach \x in {0,...,19} {\node at (1.5, \x*2.2/19 ) {$\cdot$};};
		\foreach \x in {0,...,23} {\node at (1.5, \x*2.2/23 ) {$\cdot$};};
		\foreach \x in {0,...,29} {\node at (1.5, \x*2.2/29 ) {$\cdot$};};
		
		\end{tikzpicture}
		
	\caption{As in the example, let $\Lambda_1=\mathbb{Z}+i\sqrt{2}\mathbb{Z}$ and $\Lambda_2=\mathbb{Z}+i \sqrt{5}\mathbb{Z}$. Given a point $z=x+iy \in \mathbb{C}$, the set of points of the form $(z,z+\lambda)$ for $\lambda \in \Lambda_2$ is dense in the set of points $(z, \Re(z)+y)$ for $y \in \mathbb{R}$. This is shown in this figure: for a fixed $z \in \mathbb{C}/\Lambda_1$ on the left, we see infinitely many determinations on the right, filling the vertical line.}\label{abeliandetermination}

	\end{center}
\end{figure}
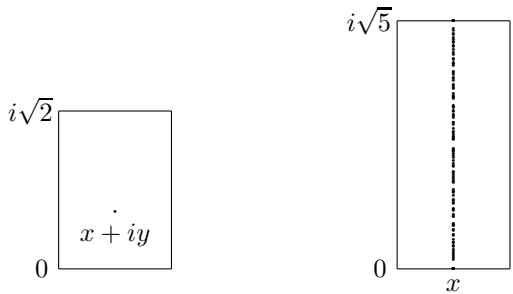

The homology argument then implies that $\o{\exp_{E_1 \times E_2}(L)} \cap W \neq \varnothing$: in algebraic terms, there are $x+iy_1$ and $x+iy_2$ in $\mathbb{C}$ such that $F(\exp_{E_1 \times E_2}(x+iy_1,x+iy_2))=0$, simply because $\o{\exp_{E_1 \times E_2}(L)}$ and $W$ are closed subsets of $E_1 \times E_2$ and the intersection product of their homology classes is non-zero.

Then we may conclude by density: by openness of the $\delta$-map (which in this case holds everywhere on $L \times W$, because they are both spaces of dimension and codimension 1) we find $(l,w) \in L \times W$ (with $w$ arbitrarily close to $(\exp_{E_1 \times E_2}(x+iy_1,x+iy_2))$) such that $w-\exp_{E_1 \times E_2}(l) \in \exp_{E_1 \times E_2}(L)$, which as usual implies that $w \in \exp_{E_1 \times E_2}(L)$ to begin with. Therefore, there is $z \in \mathbb{C}$ such that $w=\exp_{E_1 \times E_2}(z,z)$, giving a zero of the system of equations under discussion.

\printbibliography
\end{document}